\documentstyle[11pt,twoside]{article}
\title{Exponential asymptotics of the Voigt functions}
\author{\sc R. B.\ Paris \\
{\em School of Engineering, Computing and Applied Mathematics,} \\
{\em University of Abertay Dundee, Dundee DD1 1HG, UK}}
\begin{document}
\def\f#1#2{\mbox{${\textstyle \frac{#1}{#2}}$}}
\def\dfrac#1#2{\displaystyle{\frac{#1}{#2}}}
\def\boldal{\mbox{\boldmath $\alpha$}}
\newcommand{\bee}{\begin{equation}}
\newcommand{\ee}{\end{equation}}
\newcommand{\lam}{\lambda}
\newcommand{\ka}{\kappa}
\newcommand{\al}{\alpha}
\newcommand{\th}{\theta}
\newcommand{\fr}{\frac{1}{2}}
\newcommand{\fs}{\f{1}{2}}
\newcommand{\g}{\Gamma}
\newcommand{\br}{\biggr}
\newcommand{\bl}{\biggl}
\newcommand{\ra}{\rightarrow}
\newcommand{\gtwid}{\raisebox{-.8ex}{\mbox{$\stackrel{\textstyle >}{\sim}$}}}
\newcommand{\ltwid}{\raisebox{-.8ex}{\mbox{$\stackrel{\textstyle <}{\sim}$}}}
\renewcommand{\topfraction}{0.9}
\renewcommand{\bottomfraction}{0.9}
\renewcommand{\textfraction}{0.05}
\newcommand{\mcol}{\multicolumn}
\date{}
\maketitle
\pagestyle{myheadings}
\markboth{\hfill \sc R. B.\ Paris  \hfill}
{\hfill \sc Asymptotics of the Voigt functions\hfill}
\begin{abstract}
We obtain the asymptotic expansion of the Voigt functions $K(x,y)$ and $L(x,y)$ for large (real) values of the variables $x$ and $y$, paying particular attention to the exponentially small contributions. A Stokes phenomenon is encountered as $y\rightarrow+\infty$ with $x>0$ fixed.
Numerical examples are presented to demonstrate the accuracy of these new expansions.

\vspace{0.3cm}
\noindent {\bf Mathematics subject classification (2010):} 30E15, 33E20, 33C15, 34E05, 41A60
\vspace{0.1cm}
 
\noindent {\bf Keywords:} Voigt functions, asymptotic expansions, exponential asymptotics
\end{abstract}

\vspace{0.3cm}

\noindent $\,$\hrulefill $\,$

\vspace{0.2cm}

\begin{center}
{\bf 1. \  Introduction}
\end{center}
\setcounter{section}{1}
\setcounter{equation}{0}
\renewcommand{\theequation}{\arabic{section}.\arabic{equation}}
Several different notations and normalisations for the Voigt functions exist in the literature; here we adopt the 
notation employed in \cite{HJ, PS} and define the Voigt functions $K(x,y)$ and $L(x,y)$ by the convolution integrals
\bee\label{e11}
K(x,y)=\frac{y}{\pi}\int_{-\infty}^\infty \frac{e^{-t^2}}{(x-t)^2+y^2}dt=\frac{1}{\pi}\int_{-\infty}^\infty e^{-(x-yt)^2}\,\frac{dt}{1+t^2}
\ee
and 
\bee\label{e12}
L(x,y)=\frac{1}{\pi}\int_{-\infty}^\infty \frac{(x-t) e^{-t^2}}{(x-t)^2+y^2}dt=\frac{1}{\pi}\int_{-\infty}^\infty e^{-(x-yt)^2}\,\frac{t dt}{1+t^2},
\ee
where $x$ and $y$ are real variables.
These functions play an important role in several fields of physics, such as astrophysical spectroscopy and Doppler-broadened absorption resonances in neutron physics. They also describe
the one-dimensional decay in an infinite viscous fluid of an initial laminar velocity distribution in the $x, z$-plane of the form $1/(1+y^2)$ and $y/(1+y^2)$ \cite[p. 621]{L} and the distribution of temperature in one-dimensional heat conduction with similar initial temperature profiles \cite[p. 566]{JJ}.

Alternative integral representations for $K(x,y)$ and $L(x,y)$ are given by
\[K(x,y)=\frac{1}{\surd\pi}\int_0^\infty e^{-yt-\frac{1}{4}t^2} \cos (xt)\,dt,\quad
L(x,y)=\frac{1}{\surd\pi}\int_0^\infty e^{-yt-\frac{1}{4}t^2} \sin (xt)\,dt\]
when $y>0$, and
\[K(x,y)=\frac{y e^{-x^2}}{\surd\pi}\int_0^\infty \frac{F(t)}{(1+t)^{1/2}}dt,\quad
L(x,y)=\frac{x e^{-x^2}}{\surd\pi}\int_0^\infty \frac{F(t)}{(1+t)^{3/2}}dt,\]
where $F(t):=\exp [-y^2t+x^2/(1+t)]$.
The first pair of integrals shows that $K(x,y)$ and $L(x,y)$ are Fourier sine and cosine transforms
and the second pair of integrals was derived in \cite{HJ}; see also \cite{P}. Integrals of the Mellin-Barnes type involving gamma functions in the integrand have been given in \cite{PS}.
By extension of the Voigt functions to complex values of $x$ and $y$, several doubly infinite integrals have been evaluated in \cite{KC, R}. A list of over 80 integrals of the Voigt functions with respect to both variables is tabulated in \cite{P}.
 
Our aim in this paper is to obtain exponentially accurate asymptotics of the Voigt functions $K(x,y)$ and $L(x,y)$ as either $x$ or $y$, or both, become large.
The interest in exponentially precise asymptotics during the past two decades has shown that retention of exponentially small terms, previously neglected in asymptotics, is essential for a high-precision description. An early example that illustrated the advantage of retaining exponentially small terms in the asymptotic expansion of a certain integral was given in
\cite[p.~76]{Ob}. Although such terms are negligible in the Poincar\'e sense, their inclusion can significantly improve the numerical accuracy. In addition, the determination of exponentially small terms in the asymptotics of the Voigt functions may have a physical interpretation in applications.

To achieve this exponential improvement, we must optimally truncate
the standard algebraic expansions for $K(x,y)$ and $L(x,y)$ at, or near, their least terms so that the accuracy thus attained is comparable with the exponentially small terms; for a discussion of this process, see \cite[Chapter 6]{PK}. The remainder in the expansion is expressed exactly in terms of a single, so-called {\it terminant function}. The known asymptotics of this last function can then be exploited to determine the exponentially small contributions contained in the Voigt functions. 
We shall find that the behaviour of this exponentially small contribution is different in the limits $x\rightarrow\infty$, $y$ finite and $y\rightarrow\infty$, $x$ finite, since the former limit is associated with a Stokes phenomenon whereas the second limit is not.
In the final section we present some examples to illustrate the accuracy of the expansions obtained.
\vspace{0.6cm}

\begin{center}
{\bf 2. \ The algebraic asymptotic expansions} 
\end{center}
\setcounter{section}{2}
\setcounter{equation}{0}
\renewcommand{\theequation}{\arabic{section}.\arabic{equation}}
From the definition of the Voigt functions $K(x,y)$ and $L(x,y)$ in Section 1, it follows that
\[K(-x,y)=-K(x,-y)=K(x,y),\qquad L(-x,y)=-L(x,-y)=-L(x,y),\]
so that it is sufficient to consider only values of the variables satisfying $x\geq 0$, $y\geq 0$.
It is well known that $K(x,y)$ and $L(x,y)$ can be represented as the real and imaginary parts of simple
special functions of a complex variable, viz.
\begin{equation}\label{e21}
K(x,y)-iL(x,y)=e^{w^2} \mbox{erfc}\,w=\frac{e^{w^2}}{\surd\pi} \Gamma(\f{1}{2},w^2), \qquad w:=y+ix,
\end{equation}
where $\mbox{erfc}\,z$ and $\Gamma(a,z)$ denote
respectively the complementary error function and the incomplete gamma function \cite[pp.~160, 177]{DLMF}.
The Voigt functions can also be related to the confluent hypergeometric function ${}_1F_1$ by
\bee\label{e21a}
K(x,y)-iL(x,y)=e^{w^2}-\frac{2w}{\surd\pi}\,{}_1F_1(1;\f{3}{2};w^2).
\ee

From this last expression we obtain the special values:
\[K(x,0)=e^{-x^2},\qquad K(0,y)=e^{y^2} \mbox{erfc}\,y,\]
\bee\label{e21b}
L(x,0)=\mbox{Im}\,[e^{-x^2} \mbox{erf}\,(ix)]=\frac{2xe^{-x^2}}{\surd\pi}\,{}_1F_1(\fs; \f{3}{2}; x^2),\qquad L(0,y)=0.
\ee
It can also be seen from (\ref{e21a}) that the asymptotics of $K(x,y)$ and $L(x,y)$ will be different in the limits $y\ra\infty$, $x$ finite and $x\ra\infty$, $y$ finite (when $\arg\,w\simeq 0$
and $\arg\,w\simeq\fs\pi$, respectively), since the ray $\arg\,z=\pi$ is a Stokes line for the  function ${}_1F_1(a;b;z)$. The exponentially small expansions of the confluent hypergeometric functions on the negative real axis have been recently discussed in \cite{CHF}.

Although the asymptotic expansions of $K(x,y)$ and $L(x,y)$ can be obtained from the exponentially improved asymptotics of the confluent hypergeometric function ${}_1F_1(a;b;z)$ \cite[(13.2.41), (13.7.10)]{DLMF}, we find it more straightforward to employ repeated application of the recurrence relation 
\[\g(a+1,z)=a \g(a,z)+z^a e^{-z}\]
to the incomplete gamma function in (\ref{e21}). For arbitrary positive integer $m$, we obtain
\[\frac{1}{\surd\pi} \g(\fs,w^2)=
\frac{e^{-w^2}}{\surd\pi}\sum_{k=0}^{m-1} \frac{(-)^k (\fs)_k}{w^{2k+1}}+Q(\fs-m,w^2),\]
where $(a)_k=\g(a+k)/\g(a)$ is Pochhammer's symbol and $Q(a,z)=\g(a,z)/\g(a)$ is the normalised incomplete gamma function; compare \cite[(8.4.14)]{DLMF}.
If we now introduce the so-called {\it terminant function} $T_\nu(z)$ defined\footnote{In \cite{O} a different normalisation was employed for $T_\nu(z)$.} by \cite{O}
(see also \cite[p.~260]{PK} where it is called ${\hat T}_\nu(z)$)
\bee\label{e22}
T_\nu(z)=e^{\pi i\nu} \frac{\g(\nu)}{2\pi i}\,\g(1-\nu,z),
\ee
then some straightforward algebra shows that $Q(\fs-m,w^2)=T_{m+\fr}(w^2)$. Hence we obtain the result
\bee\label{e23}
K(x,y)-iL(x,y)=\frac{1}{\surd\pi}\sum_{k=0}^{m-1} \frac{(-)^k (\fs)_k}{w^{2k+1}}+2e^{w^2}T_{m+\fr}(w^2)
\ee
for arbitrary positive integer $m$. 

It should be emphasised that (\ref{e23}) is {\it exact\/} and that no approximation has yet been involved.
The summation index $m$ in the series on the right-hand side of (\ref{e23}) will now be chosen to correspond to optimal truncation --- that is, truncation at or near the term of least magnitude.
This is easily verified to be $m\simeq |w|^2$. With $m$ so chosen, Olver \cite{O} has shown that \[e^{w^2}T_{m+\fr}(w^2)=O(e^{-w^2-|w|^2}) \qquad (|w|\ra\infty,\ \ |\arg\,w|\leq\fs\pi),\]
so that in the sector $0\leq\arg\,w\leq\fs\pi-\delta$, where throughout $\delta$ denotes an abitrarily small positive quantity, the second term on the right-hand side of (\ref{e23}) is exponentially small.
Consequently, we obtain the known algebraic asymptotic expansions  given by 
\begin{eqnarray}
K(x,y)&\sim&\frac{1}{\surd\pi}\sum_{k=0}^\infty\frac{(-)^k(\fs)_k \cos (2k+1)\theta}{(x^2+y^2)^{k+\fr}}, \label{e23a}\\
L(x,y)&\sim& \frac{1}{\surd\pi}\sum_{k=0}^\infty\frac{(-)^k(\fs)_k \sin (2k+1)\theta}{(x^2+y^2)^{k+\fr}}\label{e23b}
\end{eqnarray}
as $|w|\ra\infty$ in $|\arg\,w|\leq\fs\pi-\delta$, where
\[\theta=\arg\,w=\arctan (x/y).\]
The expansions (\ref{e23a}) and (\ref{e23b}) hold when either $x$ or $y$, or both, become large in the stated sector. 
These expansions have been obtained in \cite{DA} and \cite{P} using different methods; the first two terms were derived in \cite{HJ} from an integral representation.

\vspace{0.6cm}

\begin{center}
{\bf 3. \ Exponentially improved expansions}
\end{center}
\setcounter{section}{3}
\setcounter{equation}{0}
\renewcommand{\theequation}{\arabic{section}.\arabic{equation}}
To determine the exponentially improved expansions of the Voigt functions we require the expansion of $T_\nu(z)$ as $|z|\ra\infty$ when 
$\nu=|z|+\alpha$, with $|\alpha|$ bounded. If we let $\psi=\arg\,z$, we have the integral representation for $T_\nu(z)$ in the form \cite{O}
\bee\label{e25}
e^zT_\nu(z)=\frac{e^{(\pi-\psi)i\nu}}{2\pi i} \int_0^\infty e^{-|z|\tau}\,\frac{\tau^{\nu-1}}{1+\tau e^{-i\psi}}\,d\tau
\ee
valid when $|\arg\,z|<\pi$. The asymptotic behaviour of the above integral as $|z|\ra\infty$ is governed by a saddle point at $\tau=1$ which becomes coincident with the pole of the integrand at $\tau=-e^{i\psi}$ on the Stokes lines $\psi=\pm\pi$. Olver \cite{O} has shown that for large $|z|$ with $\nu\simeq |z|$ 
\begin{eqnarray}
T_\nu(z)&\sim& -\frac{ie^{i\phi\nu}}{1-e^{i\phi}}\,\frac{e^{-z-|z|}}{\sqrt{2\pi |z|}}\sum_{k=0}^\infty
A_{2k}(\phi,\alpha) |z|^{-k},\qquad (|\arg\,z|\leq\pi-\delta)\label{e26a}\\
T_\nu(z)&\sim&\frac{1}{2}\mbox{erfc}\,[c(\phi)(\fs|z|)^\fr]-\frac{ie^{-z-|z|+i\phi|z|}}{\sqrt{2\pi |z|}}
\sum_{k=0}^\infty B_{2k}(\phi,\alpha) |z|^{-k} \nonumber\\
&& \hspace{5cm}(-\pi+\delta\leq\arg\,z\leq 3\pi-\delta),\label{e26b}
\end{eqnarray}
where $\phi=\pi-\psi$. 

The quantity $c(\phi)$, which measures the proximity of the saddle point $\tau=1$ to the pole at $\tau=e^{i\phi}$, is defined implicitly by
\bee\label{e27}
\fs c^2(\phi)=1-i\phi-e^{-i\phi}
\ee
and corresponds to the branch of $c(\phi)$ that behaves near $\phi=0$ ($\psi=\pi$) like
\[c(\phi)=\phi-\f{1}{6}i\phi^2-\f{1}{36}\phi^3+\cdots\ .\]
The coefficients $A_{2k}(\phi,\alpha)$ and $B_{2k}(\phi,\alpha)$ are given by
\[A_0(\phi,\alpha)=1,\qquad A_2(\phi,\alpha)=\f{1}{12}+h_2(\phi,\alpha),\]
\bee\label{e27a}
A_4(\phi,\alpha)=\f{1}{288}+\f{1}{12}h_2(\phi,\alpha)+2h_3(\phi,\alpha)+3h_4(\phi,\alpha), \ldots \ ,
\ee
and
\bee\label{e27b}
B_{2k}(\phi,\alpha)=\frac{e^{i\phi\alpha}A_{2k}(\phi,\alpha)}{1-e^{i\phi}}-\frac{i(-)^k 2^k (\fs)_k}{(c(\phi))^{2k+1}},
\ee
where $h_k(\phi,\alpha)$ is defined in (\ref{a4a}). Higher coefficients can be obtained from (\ref{a4}) and Table 3 in the Appendix.
With the help of {\it Mathematica} the values of $B_{2k}(\phi,\alpha)$ (for $k\leq 2$) when $\phi\ra 0$ are found to be
\[B_0(\phi,\alpha)=\f{2}{3}-\alpha-\f{1}{12}i(1-6\alpha+6\alpha^2)\phi+O(\phi^2),\]
\bee\label{e27c}
B_2(\phi,\alpha)=\f{23}{270}-\f{5}{12}\alpha+\fs a^2-\f{1}{6}a^3+O(i\phi),
\ee
\[B_4(\phi,\alpha)=\f{23}{3024}-\f{21}{160}\alpha+\f{3}{8}\alpha^2-\f{7}{18}\alpha^3+\f{1}{6}\alpha^4-\f{1}{40}\alpha^5+O(i\phi).\]

The expansion (\ref{e26a}) is obtained from a straightforward application of Laplace's method to the integral (\ref{e25}) as discussed in the Appendix, with the expansion (\ref{e26b}) resulting from the use of a standard quadratic change of variables in (\ref{e25}) to deal with the saddle and pole of the integrand; for details, see \cite{O} and \cite[\S 6.2.6]{PK}. 
In the sector $|\arg\,z|\leq\pi-\delta$, the quantity $c(\phi)$ lies in the fourth quadrant bounded away from the origin, so that the expansion of erfc\,$z$ given by \cite[p.~164]{DLMF}
\bee\label{e27d}
\mbox{erfc}\,z\sim\frac{e^{-z^2}}{\surd\pi} \sum_{k=0}^\infty (-)^k (\fs)_kz^{-2k-1}\qquad (|z|\ra\infty,\ \ |\arg\,z|<\f{3}{4}\pi)
\ee
can be used to show that the expansion (\ref{e26b}) reduces to (\ref{e26a}) away from a neighbourhood of $|\arg\,z|=\pm\pi$. On the Stokes line $\arg\,z=\pi$ ($\phi=0$), it is seen that $T_\nu(z)=\fs+O(z^{-\fr})$. In the region $\pi+\delta\leq\arg\,z\leq\pi-\delta$ enclosing the Stokes line the smooth transition of the leading behaviour of $T_\nu(z)$ (with $\nu\simeq |z|$) is described by the error function.

If we now set the parameter $\nu$ appearing in the terminant function in (\ref{e23}) equal to
\bee\label{e300}
\nu=m+\fs=|w|^2+\alpha,
\ee
where $|\alpha|$ is bounded, we find from (\ref{e26a}) the expansion
\[2e^{w^2} T_{m+\fr}(w^2)\sim\frac{e^{-x^2-y^2}}{\sqrt{2\pi} \cos \theta}\,\sum_{k=0}^\infty \frac{e^{im\phi} A_{2k}(\phi,\alpha)}{(x^2+y^2)^{k+\fr}}\qquad (\phi=\pi-2\theta)\]
as $|w|\ra\infty$ in $|\arg\,w|\leq\fs\pi-\delta$.
Substitution of this last result into (\ref{e23}) then yields the following theorem.
\newtheorem{theorem}{Theorem}
\begin{theorem} Let $w=x+iy$ be a complex variable situated in the first quadrant, 
with 
$\theta=\arctan (x/y)$ and $\phi=\pi-2\theta$. Then, if the truncation index $m$ is chosen to satisfy $m=|w|^2+\alpha-\fs$,
with $|\alpha|$ bounded, we have the compound expansions
\bee\label{e31a}
K(x,y)\sim\frac{1}{\surd\pi}\sum_{k=0}^{m-1}\frac{(-)^k(\fs)_k \cos (2k+1)\theta}{(x^2+y^2)^{k+\fr}}
+\frac{e^{-x^2-y^2}}{\sqrt{2\pi}\,\cos \theta}\sum_{k=0}^\infty \frac{\mbox{Re}\,[e^{im\phi} A_{2k}(\phi, \alpha)]}{(x^2+y^2)^{k+\fr}},\hspace{0.7cm}
\ee
\bee\label{e31b}
L(x,y)\sim\frac{1}{\surd\pi}\sum_{k=0}^{m-1}\frac{(-)^k(\fs)_k \sin (2k+1)\theta}{(x^2+y^2)^{k+\fr}}
-\frac{e^{-x^2-y^2}}{\sqrt{2\pi}\,\cos \theta}\sum_{k=0}^\infty \frac{\mbox{Im}\,[e^{im\phi} A_{2k}(\phi, \alpha)]}{(x^2+y^2)^{k+\fr}}
\ee
as $|w|\ra\infty$ in the sector $0\leq\arg\,w\leq\fs\pi-\delta$. The first few coefficients $A_{2k}(\phi,\alpha)$ are defined in (\ref{e27a}).
\end{theorem}

The expansions (\ref{e31a}) and (\ref{e31b}) break down as $\phi\ra 0$ ($\theta\ra\fs\pi$) due to the 
unbounded nature of the coefficients $A_{2k}(\phi,\alpha)$, which is caused by the presence of the Stokes line 
on $\arg\,w=\fs\pi$ of the complementary error function in (\ref{e21}). 
From (\ref{e26b}) we obtain the expansion (when $m\simeq |w|^2$)
\[2e^{w^2} T_{m+\fr}(w^2)\sim \frac{e^{-x^2-y^2}}{\sqrt{2\pi}}\left\{E(\phi)-2i\sum_{k=0}^\infty \frac{B_{2k}(\phi,\alpha)}{(x^2+y^2)^{k+\fr}}\right\}e^{i(m+\fr-\alpha)\phi}\]
for $|w|\ra\infty$ in the sector $0\leq\arg\,w\leq\fs\pi$, where we have defined
\bee\label{e35}
E(\phi):=\sqrt{2\pi}\, e^{\zeta^2} \mbox{erfc}\,\zeta, \qquad \zeta:=c(\phi) |w|/\surd 2,
\ee
which has the limiting value $E(0)=(2\pi)^{1/2}$. 
The expansion of the Voigt functions valid in a neighbourhood of $\arg\,w=\fs\pi$ is obtained by substitution of this last result into (\ref{e23}) to yield the following theorem.
\begin{theorem} 
Let $w=x+iy$ be a complex variable situated in the first quadrant, 
with $\theta=\arctan (x/y)$ and $\phi=\pi-2\theta$. Let $E(\phi)$ be as defined in (\ref{e35}). Then, if the truncation index $m$ is chosen to satisfy $m=|w|^2+\alpha-\fs$,
with $|\alpha|$ bounded, we have the compound expansions
\[K(x,y)\sim\frac{1}{\surd\pi}\sum_{k=0}^{m-1}\frac{(-)^k(\fs)_k \cos (2k+1)\theta}{(x^2+y^2)^{k+\fr}}\hspace{5cm}\]
\bee\label{e32a}
\hspace{3cm}+\frac{e^{-x^2-y^2}}{\sqrt{2\pi}}\mbox{Re}\left\{e^{i(m+\fr-\alpha)\phi} E(\phi) + \sum_{k=0}^\infty \frac{e^{im\phi}{\hat B}_{2k}(\phi,\alpha)}{(x^2+y^2)^{k+\fr}}\right\},
\ee
\[L(x,y)\sim\frac{1}{\surd\pi}\sum_{k=0}^{m-1}\frac{(-)^k(\fs)_k \sin (2k+1)\theta}{(x^2+y^2)^{k+\fr}}\hspace{5cm}\]
\bee\label{e32b}
\hspace{3cm}-\frac{e^{-x^2-y^2}}{\sqrt{2\pi}}\mbox{Im}\left\{e^{i(m+\fr-\alpha)\phi} E(\phi) + \sum_{k=0}^\infty \frac{e^{im\phi}{\hat B}_{2k}(\phi,\alpha)}{(x^2+y^2)^{k+\fr}}\right\},
\ee
as $|w|\ra\infty$ in the sector $0\leq\arg\,w\leq\fs\pi$, where 
\bee\label{e3c}
{\hat B}_{2k}(\phi,\alpha)=-2i e^{i\phi(\fr-\alpha)} \,B_{2k}(\phi,\alpha)=
\frac{A_{2k}(\phi,\alpha)}{\cos \theta}\,-\frac{(-)^k 2^{k+1} (\fs)_k}{(c(\phi))^{2k+1}}\,e^{i\phi(\fr-\alpha)}
\ee
and the first few coefficients $A_{2k}(\phi,\alpha)$ are given in (\ref{e27a}).  The quantity $c(\phi)$ is defined
in (\ref{e27}) with the branch chosen such that $c(\phi)\sim\phi$ as $\phi\ra 0$.
\end{theorem}
As discussed in Section 2, the expansions (\ref{e32a}) and (\ref{e32b}) go over into those in (\ref{e31a}) and (\ref{e31b}) away from the neighbourhood of the Stokes line $\arg\,w=\fs\pi$. 

\vspace{0.6cm}

\begin{center}
{\bf 4. \ Numerical results and concluding remarks}
\end{center}
\setcounter{section}{4}
\setcounter{equation}{0}
\renewcommand{\theequation}{\arabic{section}.\arabic{equation}}
We present in this section some numerical results to illustrate the accuracy of the expansions in Theorems 1 and 2.
We define the quantities
\[{\hat K}(x,y):=K(x,y)-\frac{1}{\surd\pi}\sum_{k=0}^{m-1}\frac{(-)^k(\fs)_k \cos (2k+1)\theta}{(x^2+y^2)^{k+\fr}}\]
and
\[{\hat L}(x,y):=L(x,y)-\frac{1}{\surd\pi}\sum_{k=0}^{m-1}\frac{(-)^k(\fs)_k \sin (2k+1)\theta}{(x^2+y^2)^{k+\fr}}.\]
Then from Theorem 1 we obtain the expansions
\bee\label{e41}
\left.\begin{array}{r}{\hat K}(x,y)\\{\hat L}(x,y)\end{array}\!\!\right\}\sim
\pm\frac{e^{-x^2-y^2}}{\sqrt{2\pi}\,\cos \theta}\left\{\!\!\!\begin{array}{c}\mbox{{Re}}\\ \mbox{{Im}}\end{array}\!\!\!\right\}\, \sum_{k=0}^\infty \frac{
e^{im\phi}A_{2k}(\phi,\alpha)}{(x^2+y^2)^{k+\fr}}
\ee
as $|w|\ra\infty$ in the sector $0\leq\arg\,w\leq\fs\pi-\delta$,
where $m$ is the optimal truncation index defined in (\ref{e300}). Similarly, from Theorem 2, we obtain the expansions
\bee\label{e42}
\left.\begin{array}{r}{\hat K}(x,y)\\{\hat L}(x,y)\end{array}\!\!\right\}\sim
\pm\frac{e^{-x^2-y^2}}{\sqrt{2\pi}}\left\{\!\!\!\begin{array}{c}\mbox{{Re}}\\ \mbox{{Im}}\end{array}\!\!\!\right\}
\left[e^{i(m+\fr-\alpha)\phi}E(\phi)+\sum_{k=0}^\infty \frac{
e^{im\phi}{\hat B}_{2k}(\phi,\alpha)}{(x^2+y^2)^{k+\fr}}\right]
\ee
as $|w|\ra\infty$ in the sector $0\leq\arg\,w\leq\fs\pi$, where $E(\phi)$ is given in (\ref{e35}).
\begin{table}[bh]
\begin{center}
\caption{\footnotesize{Values of ${\hat K}(x,y)$ and ${\hat L}(x,y)$ using the asymptotic expansions (\ref{e41}) and (\ref{e42}) for different truncation index $k$ and two values of $\theta=\arg\,w$ when $|w|=|y+ix|=3.5$. The optimal truncation index is $m=12$ ($\alpha=0.25$). The exact value is presented at the foot of each column.}}
\begin{tabular}{c|c|c||c|c}
\mcol{1}{c|}{} & \mcol{2}{c||}{$\theta=\pi/10$,\ Eq.\ (4.1)} & \mcol{2}{c}{$\theta=3\pi/8$,\ Eq.\ (4.2)}\\
\mcol{1}{c|}{$k$} & \mcol{1}{c|}{${\hat K}(x,y)$} &\mcol{1}{c||}{${\hat L}(x,y)$} & \mcol{1}{c|}{${\hat K}(x,y)$} &\mcol{1}{c}{${\hat L}(x,y)$}\\
[.1cm]\hline
&&&&\\[-0.25cm]
0 & $+1.77219153(-7)$ & $+5.45424470(-7)$ & $-1.30341265(-6)$ & $-7.12131744(-8)$\\
1 & $+1.73069912(-7)$ & $+5.51113801(-7)$ & $-1.30412205(-6)$ & $-7.19080750(-8)$\\
2 & $+1.73151197(-7)$ & $+5.50673634(-7)$ & $-1.30410921(-6)$ & $-7.18533969(-8)$\\
3 & $+1.73163893(-7)$ & $+5.50694322(-7)$ & $-1.30410846(-6)$ & $-7.18527877(-8)$\\
4 & $+1.73161147(-7)$ & $+5.50694282(-7)$ & $-1.30410848(-6)$ & $-7.18528623(-8)$\\
[.2cm]\hline
&&\\[-0.25cm]
 & $+1.73161445(-7)$ & $+5.50694067(-7)$ & $-1.30410848(-6)$ & $-7.18528635(-8)$\\
\end{tabular}
\end{center}
\end{table}

In Table 1 we show the values\footnote{In Tables 1 and 2 we write the values as $x(y)$ instead of $x\times 10^y$.} obtained from the right-hand sides of (\ref{e41}) and (\ref{e42}) for different truncation index $k$ and two values of $\theta=\arg\,w$ when $|w|=3.5$. The corresponding exact values of ${\hat K}(x,y)$ and ${\hat L}(x,y)$ are shown at the foot of each column. Table 2 shows the absolute relative error in the computation of ${\hat K}(x,y)$ and ${\hat L}(x,y)$ for a fixed $|w|$ and different argument $\theta$ using a truncation index $k=2$ in the expansions
(\ref{e41}) and (\ref{e42}). It is apparent that (\ref{e41}) breaks down as $\theta\ra\fs\pi$ whereas (\ref{e42}) holds uniformly throughout the range $0\leq\theta\leq\fs\pi$. 
\begin{table}[t]
\begin{center}
\caption{\footnotesize{The absolute relative error in the computation of ${\hat K}(x,y)$ and ${\hat L}(x,y)$ using the asymptotic expansions (\ref{e41}) and (\ref{e42}) for different $\theta=\arg\,w$ when $|w|=|y+ix|=6$. The truncation index used in the exponentially small expansions in each case is $k=2$ and $m=36$ ($\alpha=0.50$).}}
\begin{tabular}{l|c|c||c|c}
\mcol{1}{c|}{} & \mcol{2}{c||}{Eq. (4.1)} & \mcol{2}{c}{Eq. (4.2)}\\
\mcol{1}{c|}{$\theta/\pi$} & \mcol{1}{c|}{${\hat K}(x,y)$} &\mcol{1}{c||}{${\hat L}(x,y)$}
& \mcol{1}{c|}{${\hat K}(x,y)$} &\mcol{1}{c}{${\hat L}(x,y)$}\\
[.1cm]\hline
&&&&\\[-0.25cm]
0 & 1.107$(-6)$ & $-$ & 5.785$(-7)$ & $-$\\
0.10 & 1.331$(-6)$ & 2.861$(-6)$ & 4.063$(-8)$ & 4.523$(-7)$\\
0.20 & 1.782$(-5)$ & 6.846$(-6)$ & 3.325$(-7)$ & 9.146$(-8)$\\
0.30 & 2.274$(-4)$ & 3.882$(-5)$ & 2.082$(-7)$ & 9.980$(-9)$\\
0.40 & 1.282$(-3)$ & 6.554$(-3)$ & 2.215$(-8)$ & 3.897$(-8)$\\
0.45 & 3.111$(-1)$ & 1.298$(-1)$ & 5.132$(-8)$ & 5.343$(-9)$\\
0.48 & 5.120$(+0)$ & 1.434$(+1)$ & 4.116$(-8)$ & 1.647$(-9)$\\
[.2cm]\hline
\end{tabular}
\end{center}
\end{table}

From (\ref{e41}) and the fact that $A_0(\phi,\alpha)=1$, the leading behaviour of the exponentially small contribution (when the algebraic expansions are optimally truncated according to (\ref{e300})) is
\[\left.\begin{array}{r}{\hat K}(x,y)\\{\hat L}(x,y)\end{array}\!\!\right\}\sim\pm\frac{e^{-x^2-y^2}\sec \theta}{\sqrt{2\pi(x^2+y^2)}}\,\left\{\!\!\begin{array}{c}\cos m\phi\\ \sin m\phi\end{array}\!\!\right\}
=\frac{(-)^m e^{-x^2-y^2}}{\sqrt{2\pi}\,y}\,\left\{\!\!\begin{array}{c}\cos 2m\theta\\ \sin 2m\theta\end{array}\!\!\right\}\]
for $|w|\ra\infty$ away from the Stokes line $\arg\,w=\fs\pi$.
From (\ref{e42}) and the small-$\phi$ behaviour of $B_0(\phi,\alpha)$ in (\ref{e27c}), the leading behaviour in the neighbourhood of $\arg\,w=\fs\pi$ ($\phi\simeq 0$) is found to be approximately
\[\left.\begin{array}{r}{\hat K}(x,y)\\{\hat L}(x,y)\end{array}\!\!\right\}\sim\frac{e^{-x^2-y^2}}{\sqrt{2\pi}}
\left[\pm\left\{\!\!\!\begin{array}{c}\mbox{{Re}}\\ \mbox{{Im}}\end{array}\!\!\!\right\}[e^{i(m+\fr-\alpha)\phi}E(\phi)] \right.\]
\[\left.+\frac{1}{\sqrt{x^2+y^2}}\left\{\bl(\frac{4}{3}-2\alpha\br)\left\{\!\!\begin{array}{c}\sin m\phi\\ \cos m\phi\end{array}\!\!\right\}
\pm\bl(\frac{1}{2}-\frac{4}{3}\alpha+\alpha^2\br)\left\{\!\!\begin{array}{c}\phi \,\cos m\phi\\ \phi\, \sin m\phi\end{array}\!\!\right\}+O(\phi^2)\right\}\right].\]

Finally, when $x=0$, $y>0$ ($\theta=0$, $\phi=\pi$), the coefficients $A_{2k}(\pi,\alpha)$ are real and when $y=0$, $x>0$ ($\theta=\fs\pi$, $\phi=0$), we have $c(0)=0$ and the coefficients $B_{2k}(0,\alpha)$ are real\footnote{This is a conjecture. It has been verified for the coefficients with $k\leq 5$.}; see (\ref{e27c}) and the definition (\ref{e3c}). It follows from (\ref{e31a}) and (\ref{e31b}) that as $y\ra\infty$
\[K(0,y)\sim \frac{1}{\surd\pi}\sum_{k=0}^{m-1}\frac{(-)^k (\fs)_k}{y^{2k+1}}+O(y^{-1}e^{-y^2}),\qquad L(0,y)\sim 0\]
and from (\ref{e32a}) and (\ref{e32b}) that as $x\ra\infty$ 
\[K(x,0)\sim e^{-x^2},\qquad L(x,0) \sim \frac{1}{\surd\pi}\sum_{k=0}^{m-1}\frac{(\fs)_k}{x^{2k+1}}+O(x^{-1}e^{-x^2})\]
in accordance with the special values in (\ref{e21b}).
\vspace{0.6cm}

\begin{center}
{\bf Appendix: The coefficients $A_{2k}(\phi,\alpha)$}
\end{center}
\setcounter{section}{1}
\setcounter{equation}{0}
\renewcommand{\theequation}{\Alph{section}.\arabic{equation}}
We derive expressions for the coefficients $A_{2k}(\phi,\alpha)$ for $k\leq 5$ appearing in the expansion (\ref{e26a}) of the terminant function $T_\nu(z)$. This follows the procedure given in \cite{O} which we describe here for completeness. 

We consider the integral in (\ref{e25}) when $\nu=|z|+\alpha$, with $|\alpha|$ bounded, in the form
\[I=\int_0^\infty e^{-|z|(\tau-\log\,\tau)}\,\frac{\tau^{\alpha-1}}{1-\tau e^{i\phi}}\,d\tau \qquad (|\arg\,z|<\pi),\]
where $\phi=\pi-\arg\,z$. 
As $|z|\ra\infty$, the exponential factor has a saddle point at $\tau=1$. 
Accordingly, we set $\tau=1+t$ and introduce the variable
\[\fs w^2=t-\log\,(1+t), \qquad \frac{d\tau}{dw}=\frac{w\tau}{t}\]
with $w\sim t$ as $t\ra 0$,
to find upon reversion the expansions
\bee\label{a1}
t=w+\f{1}{3}w^2+\f{1}{36}w^3-\f{1}{270}w^4+ \f{1}{4320}w^5+\f{1}{17010}w^6+\cdots
\ee
and 
\bee\label{a2}
\frac{w}{t}=1-\f{1}{3}w+\f{1}{12}w^2-\f{2}{135}w^3+\f{1}{864}w^4+\f{1}{2835}w^5+\cdots.
\ee
The integral $I$ then becomes
\bee\label{a3a}
I=\int_{-\infty}^\infty e^{-\fr |z|w^2}\,\frac{\tau^{\alpha-1}}{1-\tau e^{i\phi}}\,\frac{d\tau}{dw}\,dw=
\int_{-\infty}^\infty e^{-\fr |z|w^2}\,\frac{wg(t)}{t\,}\,dw,
\ee
where
\[g(t):=\frac{\tau^\alpha}{1-\tau e^{i\phi}}.\]

The function $g(t)$ may be expanded about $t=0$ ($\tau=1$) in the form
\[g(t)=\sum_{k=0}^\infty t^k \frac{g^{(k)}(0)}{k!}=\frac{1}{1-e^{i\phi}}\sum_{k=0}^\infty t^k h_k(\phi,\alpha),
\qquad(|t|<\min \{1,2\sin \fs\phi\})\]
where 
\bee\label{a4a}
h_k(\phi,\alpha):=(1-e^{i\phi})\,\frac{g^{(k)}(0)}{k!}=\sum_{r=0}^k \left(\!\!\begin{array}{c}\alpha\\k-r\end{array}\!\!\right) \left(\frac{e^{i\phi}}{1-e^{i\phi}}\right)^r.
\ee
We can now apply Laplace's method \cite[p.~85]{Ob}, upon 
insertion of (\ref{a1}) and (\ref{a2}) into (\ref{a3a}), to obtain 

\[I\approx\frac{1}{1-e^{i\phi}}\int_{-\infty}^\infty e^{-\fr|z|w^2}\left\{1+w^2(\f{1}{12}+h_2(\phi,\alpha))\right.\hspace{3cm}\]
\[\hspace{4cm}\left.+w^4(\f{1}{864}+
\f{1}{36}h_2(\phi,\alpha)+\f{2}{3}h_3(\phi,\alpha)+h_4(\phi,\alpha))+\cdots \right\}dw,\]
where only even powers of $w$ contribute. Evaluation of the integrals then leads to the expansion
\[I\sim\frac{\sqrt{2\pi}}{1-e^{i\phi}}\sum_{k=0}^\infty A_{2k}(\phi,\alpha)\,|z|^{-k-\fr}\]
as $|z|\ra\infty$ in $|\arg\,z|\leq\pi-\delta$,
where \cite{O}
\[A_0(\phi,\alpha)=1,\qquad A_2(\phi,\alpha)=\f{1}{12}+h_2(\phi,\alpha),\]
\bee\label{a3}
A_4(\phi,\alpha)=\f{1}{288}+\f{1}{12}h_2(\phi,\alpha)+2h_3(\phi,\alpha)+3h_4(\phi,\alpha).
\ee

Continuation of this process with the help of {\it Mathematica} shows that the coefficients $A_{2k}(\phi,\alpha)$ can be expressed in the general form
\bee\label{a4}
A_{2k}(\phi,\alpha)=(-)^k \gamma_k+\sum_{j=2}^{2k} c_{j,k}\, h_j(\phi,\alpha),
\ee
where $\gamma_k$ are the Stirling coefficients \cite[p.~32]{PK}, with $\gamma_0=1$ and
\[\gamma_1=-\f{1}{12},\quad \gamma_2=\f{1}{288},\quad \gamma_3=\f{139}{51840},  \quad\gamma_4=-\f{571}{2488320},\quad\gamma_5=-\f{163879}{209018880},\ \ldots\ .\]
The $c_{j,k}$ appear as the coefficient of $w^{2k-1}$ in the expansion of $t^j$ (multiplied by the factor $2^k (\fs)_k$) and are presented in Table 3 for $k\leq 5$.
\begin{table}[t]
\begin{center}
\begin{tabular}{c|ccccccccc}
\mcol{1}{c|}{$k\backslash j$} & \mcol{1}{c}{2} & \mcol{1}{c}{3} & \mcol{1}{c}{4} & \mcol{1}{c}{5} & \mcol{1}{c}{6}
& \mcol{1}{c}{7}& \mcol{1}{c}{8}& \mcol{1}{c}{9}& \mcol{1}{c}{10}\\
[.1cm]\hline
&&&&&&&&&\\[-0.25cm]
1 & 1 & & & & &&&&\\
&&&&&&&&&\\[-0.25cm]
2 & \f{1}{12} & 2 & 3 & & &&&&\\
&&&&&&&&&\\[-0.25cm]
3 & \f{1}{288} & \f{1}{6} & \f{25}{4} & 20 & 15&&&&\\
&&&&&&&&&\\[-0.25cm]
4 & -\f{139}{51840} & \f{1}{144} & \f{49}{96} &  \f{77}{3} & \f{525}{4} & 210 & 105 && \\
&&&&&&&&&\\[-0.25cm]
5 & -\f{571}{2488320} & -\f{139}{25920} & \f{221}{17280} & \f{149}{72} & \f{12565}{96} & \f{1883}{2} & \f{9555}{4} &
2520 & 945\\
[.2cm]\hline
\end{tabular}
\end{center}
\caption{\footnotesize{The coefficients $c_{j,k}$ for $k\leq 5$.}}
\end{table}
We observe that 
\[c_{2,k}=(-)^{k-1}\gamma_{k-1},\qquad c_{3,k}=2(-)^k\gamma_{k-2}\quad \mbox{and}\quad c_{k,2k}=2^k (\fs)_k.\]

\end{document}